\documentclass{article}
\usepackage [utf8] {inputenc}
\usepackage{mathtools}
\usepackage{amsthm}
\usepackage{amssymb}

\newtheorem{theorem}{Theorem}
\newtheorem{proposition}{Proposition}

\title{A generalization of the Brauer--Fowler theorem}
\author{Saveliy V. Skresanov
\thanks{The research was carried out within the framework of the Sobolev Institute of Mathematics state contract (project FWNF-2022-0002).}}
\date{}

\begin{document}
\maketitle

\begin{abstract}
	The famous Brauer-Fowler theorem states that the order of a finite simple group
	can be bounded in terms of the order of the centralizer of an involution.
	Using the classification of finite simple groups, we generalize this theorem and prove that if a simple locally finite group has an involution
	which commutes with at most \( n \) involutions, then the group is finite and its order is bounded in terms of~\( n \) only.
	This answers a question of Strunkov from the Kourovka notebook.
\end{abstract}

\section*{Introduction}

In 1955 Brauer and Fowler proved a result which became one of the staples of the classification of finite simple groups (CFSG):
\medskip

\noindent\textbf{The Brauer-Fowler theorem}~\cite[Corollary~(2I)]{BF}
\emph{If \( G \) is a finite simple group and the centralizer of some involution has order \( n \), then \( |G| < [n(n+1)/2]! \).
In particular, there exist only a finite number of finite simple groups with a given centralizer of an involution.}
\medskip

The proof of this theorem is elementary and relies on an ingenious technique to count involutions in a group.
Since then this method was applied to obtain other results relating the structure of a group with properties
of its involutions, notable examples being the Thompson order formula and the Brauer-Wielandt formula (see~\cite[Part I, chapter 34]{GLS2}).
The Brauer-Fowler theorem itself was also widely generalized. For example, in~\cite{GR} it is proved that for a finite group \( G \)
with an involution \( t \), one has \( |G:F(G)| < |C_G(t)|^4 \), where \( F(G) \) is the Fitting subgroup of~\( G \).
Hartley proved~\cite{H} that if a nonabelian finite simple group \( G \) admits an automorphism of order \( k \) having \( n \) fixed points,
then \( |G| \) is bounded in terms of \( k \) and~\( n \).

In 1985 Strunkov showed that in some cases the order of the centralizer of an involution can be replaced by the number of involutions
in that centralizer. The following is a corollary from his results:

\begin{proposition}[{\cite[Theorem~1]{S}}]\label{strun}
	Let \( G \) be a periodic simple group, and let \( t, z \in G \) be nonconjugate involutions.
	If the sets \( C_G(t) \cap (t^G \cup z^G) \) and \( C_G(z) \cap (t^G \cup z^G) \) are finite
	and have at most \( n \) elements, then \( G \) is also finite and its order can be bounded in terms of~\( n \) only.
\end{proposition}

In particular, if \( t \) and \( z \) centralize at most \( n \) involutions each, then \( G \) is finite and its order is bounded in terms of~\( n \).

It is natural to ask if it is possible to prove a result similar to Strunkov's for one involution only (at least in the case of finite groups).
And indeed, in 11th issue of the Kourovka notebook~\cite{kourovka} Strunkov asks:
\medskip

\noindent\textbf{Question 11.96.}
\emph{Is it true that, for a given number \( n \), there exist only finitely many finite
simple groups each of which contains an involution which commutes with at most \( n \) involutions of the group?
Is it true that there are no infinite simple groups satisfying this condition?}
\medskip

For the second part of the question the answer is ``no'', since \( \mathrm{PSL}_2(\mathbb{R}) \) is an infinite simple group
which contains an involution, but does not contain a subgroup isomorphic to the Klein 4-group. Nevertheless, the answer to the
first part of the question is positive, in fact, it is true in the more general context of simple locally finite groups.

\begin{theorem}[mod CFSG]\label{main}
	Let \( G \) be a simple locally finite group, and let \( t \in G \) be an involution. If there are at most \( n \)
	involutions in \( G \) commuting with \( t \), then \( G \) is finite and its order is bounded in terms of \( n \) only.
\end{theorem}

The principal case in the proof is that of a finite simple group of Lie type. We use a result of Guralnick and Robinson~\cite{GRcomm}
to bound the rank of the group, and a combinatorial argument together with character ratio bounds to bound the field.
An explicit bound of the form \( |G| \leq C_1 \cdot n^{C_2 \cdot n^2} \) for some universal constants \( C_1, C_2 > 0 \) can be extracted from the proof,
but it is clearly far from being tight.

We end the introduction with an open problem.
Recall that by Shunkov's theorem~\cite{Sh}, a periodic group with a finite centralizer of an involution is almost solvable,
so a simple group with such an involution is finite, and hence the Brauer--Fowler theorem can be extended to periodic groups.
Proposition~\ref{strun} of Strunkov also applies to periodic groups. 
It is quite interesting whether Theorem~\ref{main} can be generalized to periodic simple groups:
\medskip

\noindent\textbf{Open problem.}
\emph{Let \( G \) be a periodic simple group, and suppose that \( G \) contains an involution which commutes only with a finite number of involutions.
Is it true that \( G \) is finite?}
\medskip

It follows from another theorem of Strunkov~\cite[Theorem~2]{S} that if an infinite 2-group \( G \) has an involution \( t \in G \)
such that \( C_G(t) \cap t^G \) is finite, then the center of \( G \) is nontrivial. In particular, such \( G \) cannot be simple and
hence a hypothetical counterexample to the above question must contain elements of odd order.

We also note that a result of Durakov~\cite[Theorem~2]{D} implies that a positive solution to this problem, even in the case when
the centralizer of an involution \( t \in G \) is a Pr\"ufer 2-group, would resolve a question of Mazurov from the Kourovka notebook~\cite[15.54]{kourovka}.
Namely, it is conjectured that in such \( G \) all elements of odd order inverted by \( t \) form a subgroup.

\section*{Proof}

First we prove the result for \emph{finite} simple groups.
Let \( G \) be a finite simple group, and let \( t \in G \) be some involution. We may assume that \( G \) is nonabelian and its order
is larger than some universal constant, in particular it is not a sporadic finite simple group.
Suppose that there are at most \( n \) involutions commuting with \( t \). We start from a result by Guralnick and Robinson.

\begin{proposition}[{\cite[Corollary 2]{GRcomm}}]
	Let \( G \) be a finite quasisimple group.
	If \( x, y \in G \) are involutions, then \( x \) commutes with some conjugate of \( y \).
\end{proposition}

It follows that there are at most \( n \) conjugacy classes of involutions in~\( G \).
If \( G \) is an alternating group of degree \( m \), then \( G \) contains at least \( \lfloor m/4 \rfloor \) conjugacy classes of involutions, and hence \( m \leq 4n+1 \),
so the order of \( G \) is bounded in terms of~\( n \).

From now on we may assume that \( G \) is a finite simple group of Lie type of rank \( r \) over the field of order~\( q \).
Note that the rank~\( r \) of \( G \) can be bounded in terms of the number of conjugacy classes of involutions (and hence in terms of~\( n \)).
Indeed, this follows from~\cite[Table~4.5.1]{GLS} for \( q \) odd, and~\cite[Sections~4--6]{AS} for \( q \) even. Now it suffices to bound \( q \) in terms of~\( n \).

Suppose that \( G \) has at least two conjugacy classes of involutions, and let \( \mathcal{K} \) be a conjugacy class of involutions not containing~\( t \).
For any \( w \in \mathcal{K} \) the group generated by \( t \) and \( w \) is a dihedral group. Its order must be divisible by \( 4 \), since otherwise
the subgroups generated by \( t \) and \( w \) would be conjugate as the Sylow subgroups of the dihedral group in question.
Hence there exists some involution \( z \) commuting with \( t \) and \( w \).
Since there are at most \( n \) involutions commuting with \( t \), by pigeonhole principle there exists an involution \( z \)
such that it commutes with at least \( |\mathcal{K}| / n \) elements from \( \mathcal{K} \), or in other words, \( |\mathcal{K} \cap C_G(z)| \geq |\mathcal{K}| / n \).
Now to deduce a bound on \( q \) we will apply what is essentially an argument about fixed-point ratios, cf.~\cite[Lemma~1.2~(iv)]{B}.

Let \( \pi : G \to \mathbb{C} \) be the permutation character of a transitive action of \( G \) on its conjugacy class \( \mathcal{K} \).
By definition of \( \pi \), for \( g \in G \) we have \( \pi(g) = |\mathcal{K} \cap C_G(g)| \). Since \( \pi(1) = |\mathcal{K}| \), we have
\[ \frac{\pi(z)}{\pi(1)} = \frac{|\mathcal{K} \cap C_G(z)|}{|\mathcal{K}|} \geq \frac{1}{n}. \]
As the action is transitive, we have a decomposition \( \pi = 1_G + \chi_1 + \dots + \chi_k \), \( k \geq 1 \), where \( 1_G \) is the trivial character
and \( \chi_1, \dots, \chi_k \) are some nontrivial irreducible ordinary characters (with repetitions allowed). We claim that for some \( i \) we have
\[ \frac{\pi(z)}{\pi(1)} \leq \frac{1 + |\chi_i(z)|}{1 + \chi_i(1)}. \]
Set \( f = \pi(z)/\pi(1) \), and assume that for all \( i = 1, \dots, k \) we have \( f(1 + \chi_i(1)) > 1 + |\chi_i(z)| \). Then
\begin{multline*}
	\pi(z) = \left| 1 + \sum_{i = 1}^k \chi_i(z) \right| \leq 1-k + \sum_{i=1}^k (1+|\chi_i(z)|) < 1-k + f \sum_{i=1}^k (1 + \chi_i(1)) =\\
	= 1-k + f(k-1 + \pi(1)) = \pi(z) - (k-1)(1-f) \leq \pi(z)
\end{multline*}
and this is a contradiction.

Therefore for some \( i \) we have \( 1/n \leq (1 + |\chi_i(z)|)/(1 + \chi_i(1)) \). Note that since \( z \) is an involution, \( \chi_i(z) \) is an integer.
If \( \chi_i(z) = 0 \), then \( \chi_i(1) \leq n-1 \) and the order of \( G \) is bounded by Jordan's theorem (one can also obtain
a bound from the main result of~\cite{LZ}). If \( \chi_i(z) \) is nonzero, then \( |\chi_i(z)| \geq 1 \) and hence \( 1/n \leq 2|\chi_i(z)|/\chi_i(1) \).
Gluck's bound~\cite{G} for character ratios gives \( |\chi_i(z)| / \chi_i(1) \leq C \cdot q^{-1/2} \) for some universal constant \( C \),
thus \( q \leq (2Cn)^2 \). Together with the bound on the rank of \( G \) this gives an upper bound on the order of \( G \) in the case when
there are at least two conjugacy classes of involutions.

Now suppose that \( G \) has only one conjugacy class of involutions. If \( q \) is even, it follows from~\cite[Theorem~2.4.1 and Table~2.4]{GLS}
that \( G \) contains an elementary abelian subgroup of order~\( q \); a root subgroup suffices in most cases. As we can always conjugate \( t \) inside
such a subgroup, it follows that \( t \) commutes with at least \( q-1 \) involutions and hence \( q \leq n+1 \) as wanted.

Finally, if \( q \) is odd and \( G \) has a unique conjugacy class of involutions, it follows from~\cite[Table~4.5.1]{GLS} that \( G \) is one of the following groups:
\[ \mathrm{PSL}_2(q),\, \mathrm{PSL}_3(q),\, \mathrm{PSL}_4(q),\, \mathrm{PSU}_3(q),\, \mathrm{PSU}_4(q),\, G_2(q),\, {}^2G_2(q),\, {}^3D_4(q). \]
If \( G = \mathrm{PSL}_2(q) \), then the centralizer of \( t \) is a dihedral group of order \( q \pm 1 \), in particular, it contains at least \( (q-1)/2 \) involutions
and \( q \leq 2n+1 \) as wanted. Now, it is sufficient to prove that \( \mathrm{PSL}_2(q) \) lies in \( G \) in all the other cases, provided
that \( q \) is large enough. Indeed, we can always conjugate an involution inside \( \mathrm{PSL}_2(q) \) and obtain the same bound on \( q \) in terms of~\( n \).

Recall that \( \Omega_3(q) \leq \mathrm{SL}_3(q) \) and \( \Omega_3(q) \simeq \mathrm{PSL}_2(q) \) is a simple group for \( q \) large enough.
Therefore \( \mathrm{PSL}_3(q) \) contains \( \mathrm{PSL}_2(q) \). The group \( \mathrm{PSL}_4(q) \) contains \( \mathrm{SL}_3(q) \) and hence \( \mathrm{PSL}_2(q) \).
A unitary group \( \mathrm{SU}_3(q) \) contains \( \Omega_3(q) \simeq \mathrm{PSL}_2(q) \), hence \( \mathrm{PSU}_3(q) \) also contains \( \mathrm{PSL}_2(q) \).
The group \( \mathrm{SU}_4(q) \) contains \( \mathrm{GU}_3(q) \), thus again \( \mathrm{PSL}_2(q) \) can be embedded into \( \mathrm{PSU}_4(q) \).

As for the exceptional groups of Lie type, \( G_2(q) \) contains \( \mathrm{SL}_3(q) \), see~\cite[Table~4.1]{W}, and hence contains \( \mathrm{PSL}_2(q) \).
The group \( ^2G_2(q) \) contains \( \mathrm{PSL}_2(q) \) by~\cite[Theorem~4.2]{W}, and \( ^3D_4(q) \) contains \( G_2(q) \) by~\cite[Theorem~4.3]{W}.
As we showed earlier, this implies an upper bound on \( q \) in terms of \( n \), hence the theorem is proved in the case when \( G \) is finite.
\medskip

Let \( f(n) \) denote the largest size of a finite simple group, which contains an involution commuting with at most \( n \) involutions.
By what we proved above, \( f(n) \) is defined correctly.

Now assume that \( G \) is an infinite simple locally finite group, and let \( t \in G \) be an involution. Let \( t_1, \dots, t_n \in G \) be the involutions
commuting with \( t \); of course, \( t \) is one of them. Let \( T \) be a finite subgroup of \( G \) containing \( t_1, \dots, t_n \).
Since \( G \) is infinite, we may choose \( T \) to be such that \( |T| > f(n) \). A simple locally finite group has a Kegel cover~\cite[Corollary~4.3]{KW},
hence there is a finite subgroup \( H \) of \( G \) and a normal subgroup \( N \) of \( H \), such that \( H/N \) is simple,
\( T \leq H \) and \( T \cap N = 1 \).

Suppose that the order of \( N \) is even. Then it contains an involution, and hence the number of involutions in \( N \) is odd.
The involution \( t \) lies in \( H \) and acts on \( N \) by conjugation, therefore there must be an involution \( z \in N \) centralized by~\( t \).
But then \( z \) is be one of \( t_1, \dots, t_n \) and thus must lie in \( T \). This is a contradiction to \( T \cap N = 1 \), therefore the order of \( N \) is odd.

By~\cite[Corollary~3.28]{Is}, since the orders of \( t \) and \( N \) are coprime, the centralizer of an involution \( tN \) in \( H/N \)
is the image of the centralizer of \( t \) in \( H \). It follows that a finite simple group \( H/N \) has an involution which
commutes with at most \( n \) involutions. By what we proved above, this implies that the order of \( H/N \) is at most \( f(n) \).
This is a contradiction, since \( f(n) < |T| \leq |H/N| \leq f(n) \). Therefore \( G \) must be finite and the theorem is proved. \qed

\section{Acknowledgements}

The author expresses his gratitude to A.S.~Mamontov and A.V.~Vasil'ev for helpful remarks improving the exposition of this paper.

\bigskip

\noindent
\emph{Saveliy V. Skresanov}

\noindent
\emph{Sobolev Institute of Mathematics, 4 Acad. Koptyug avenue, Novosibirsk, Russia}

\noindent
\emph{Email address: skresan@math.nsc.ru}

\end{document}